\documentclass[aos,MSNbibl,nameyear,dvips]{arximspdf}

\doi{10.1214/14-AOS1203} 
\volume{42}
\issue{2}
\pubyear{2014}
\firstpage{728}
\lastpage{756}

\makeatletter
\def\mid{|}

\newcommand{\rrvert}{\vert}
\newcommand{\llvert}{\vert}
\def\cal{\mathcal}
\newcommand{\eqref}[1]{(\ref{#1})}
\renewcommand{\P}{{\cal P}}
\renewcommand{\L}{{\cal L}}
\newcommand{\D}{{\cal D}}
\newproclaim{definition}{Definition}
\newproclaim{example}{Example}
\newtheorem{theorem}{Theorem}
\newtheorem{corollary}{Corollary}
\newtheorem{lemma}{Lemma}
\newtheorem{proposition}{Proposition}
\makeatother

\begin{document}
\begin{frontmatter}

\title{Dominating countably many forecasts}
\runtitle{Dominating}

\begin{aug}
\author[a]{\fnms{M.~J.} \snm{Schervish}\corref{}\ead[label=e2]{mark@cmu.edu}},
\author[b]{\fnms{Teddy} \snm{Seidenfeld}\ead[label=e1]{teddy@stat.cmu.edu}}
\and
\author[a]{\fnms{J.~B.} \snm{Kadane}\ead[label=e3]{kadane@stat.cmu.edu}}
\affiliation{Carnegie Mellon University}
\address[a]{M.~J. Schervish\\
J.~B. Kadane\\
Department of Statistics\\
Carnegie Mellon University\\
Pittsburgh, Pennsylvania 15213\\
USA\\
\printead{e2}\\
\phantom{E-mail:\ }\printead*{e3}}
\address[b]{T. Seidenfeld\\
Department of Statistics\\
\quad and Department of Philosophy\\
Carnegie Mellon University\\
Pittsburgh, Pennsylvania 15213\\
USA\\
\printead{e1}}
\runauthor{M.~J. Schervish, T. Seidenfeld and J.~B. Kadane}
\end{aug}

\received{\smonth{10} \syear{2012}}
\revised{\smonth{1} \syear{2014}}

%
\begin{abstract}
We investigate differences between a simple Dominance Principle
applied to sums of \emph{fair prices} for variables and dominance
applied to sums of \emph{forecasts} for variables scored by proper
scoring rules. In particular, we consider differences when fair prices
and forecasts correspond to finitely additive expectations and
dominance is applied with infinitely many prices and/or forecasts.
\end{abstract}

%
\begin{keyword}[class=AMS]
\kwd[Primary ]{62A01}
\kwd[; secondary ]{62C05}
\end{keyword}
\begin{keyword}
\kwd{Proper scoring rule}
\kwd{coherence}
\kwd{conglomerable probability}
\kwd{dominance}
\kwd{finitely additive probability}
\kwd{sure-loss}
\end{keyword}

\end{frontmatter}

\section{Introduction}\label{sec:intro}
The requirement that preferences are \emph{coherent} aims to make rigorous
the idea that elementary restrictions on rational preferences entail
that personal probabilities satisfy the axioms of mathematical
probability. This use of coherence as a justification of personal
probability is very well illustrated by de Finetti's (\citeyear{deFinetti1974}) approach to the
foundations of probability. De Finetti distinguished two senses
of coherence: coherence$_1$ and coherence$_2$. Coherence$_1$ requires
that probabilistic forecasts for random variables (he calls them
previsions) do not lead to
a finite set of fair contracts that, together, are uniformly dominated
by abstaining. Coherence$_2$ requires that a finite set of probabilistic
forecasts cannot be uniformly dominated under Brier (squared error)
score by a rival set of forecasts. He showed that these two senses of
coherence are equivalent in the following sense. Each version of
coherence results in using the expectation of a random variable as its
forecast. Moreover, these expectations are
based on a finitely additive probability without requiring that
personal probability is countably additive. [In Appendix \ref{app:lf},
we explain what we mean by expectations with respect to finitely
additive probabilities. These are similar in many ways, but not
identical to integrals in the sense of \citet{dunford-schwartz1958},
Chapter III.] Schervish, Seidenfeld and
Kadane (\citeyear{ssk2009b})
extended this equivalence to include a large class of strictly proper
scoring rules (not just Brier score) but for events only. The
corresponding extension to general random variables is included in the
supplemental article [Schervish, Seidenfeld and Kadane
(\citeyear{ssk2013})]. Here, we refer to the
extended sense of coherence$_2$ as coherence$_3$.

We investigate asymmetries between coherence$_1$ and coherence$_3$
reflecting differences between cases where personal probabilities are
countably additive and where personal probabilities are finitely (but
not countably) additive. We give conditions where coherence$_3$ may be
applied to assessing countably many forecasts at once, but where
coherence$_1$ cannot be applied to combining infinitely many fair
contracts. Also, we study conditional forecasts given elements of a
partition $\pi$, where the conditional forecasts are based on the
conditional probabilities
given elements of $\pi$. Each coherence criterion is violated by combining
infinitely many conditional forecasts when those
conditional forecasts are not conglomerable (see
Definition~\ref{def:nonc}) in the partition
$\pi$. Neither criterion is violated by combining infinitely many
conditional forecasts when conditional expectations satisfy the law of
total previsions (see Definition~\ref{def:ltp}) in $\pi$.

\section{Results of de Finetti}\label{sec:results} Coherence of
preference, as de Finetti [(\citeyear{deFinetti1974}), Chapter~3] formulates it, is the
criterion that a rational decision maker respects \textit{uniform
(strict) dominance}. In Section~\ref{sub:dominance}, we
explain the version of the Dominance Principle that de Finetti uses.
In Section~\ref{sub:coherence}, we review de Finetti's two versions of
coherence, with a focus on how preferences based on a finitely
additive probability are coherent.

\subsection{Dominance}\label{sub:dominance} Let $\Omega$ be a
set. The elements of $\Omega$ will be
called \emph{states} and denoted $\omega$. Random variables are
real-valued functions with domain $\Omega$, which we denote with
capital letters. Let $I$ index a set of options. Consider a
hypothetical decision problem ${\cal O}$ specified by a set of
exclusive options ${\cal O} = \{O_i\dvtx i \in I\}$. Each option $O_i$
is a random variable with the following interpretation: If $\omega$
is the state which occurs, then $O_i(\omega)$ denotes the decision
maker's loss (negative of cardinal utility) for choosing option
$O_i$. The values of $O_{i}$ (for all
$i\in I$) are defined up to a common positive affine transformation.

\begin{definition} Let $O_i$ and $O_j$
be two options from ${\cal O}$. If there exists an $\varepsilon> 0$
such that for each $\omega\in\Omega$, $O_j(\omega) > O_i(\omega) +
\varepsilon$, then option
$O_i$ \emph{uniformly strictly dominates $O_j$}.
If, for each $\omega$, $O_j(\omega) > O_i(\omega)$, we say that $O_i$
\emph{simply dominates}~$O_j$.
\end{definition}
Uniform strict dominance is clearly stricter than simple dominance.
As we explain, next, in order to permit preferences based on
maximizing finitely (and not necessarily countably) additive
expectations, de Finetti used the following Dominance Principle,
rather than
some other more familiar concepts of admissibility, for example, simple
dominance. There are additional ways to define dominance, which
we discuss further in Section~\ref{sec:summary}.

\textsc{Dominance Principle}: Let $O_i$ and $O_j$ be options in
${\cal O} $. If $O_i$ uniformly (strictly) dominates $O_j$, then
$O_j$ is an \textit{inadmissible} choice from ${\cal O}$.

\subsection{Coherence$_1$ and coherence$_2$}\label{sub:coherence}
De Finetti [(\citeyear{deFinetti1974}), Chapter~3] formulated two criteria of
\textit{coherence} that are based on the Dominance Principle.
Throughout this paper, we follow the convention of identifying events
with their indicator functions.\vspace*{-2pt}
%
\begin{definition}\label{def:coh}\label{def:condcoh}
A \emph{conditional prevision} (or \emph{conditional forecast})
$P(X|H)$ for a
random variable $X$ given a nonempty event $H$ is a fair price for
buying and selling $X$ in the
sense that, for all real $\alpha$, the option that costs the agent
$\alpha H[X-P(X|H)]$ is considered fair. [We call $P(X|\Omega)$ an
\emph{unconditional prevision} and denote it $P(X)$.]
A collection $\{P(X_i|H_i)\dvtx i\in
I\}$ of such conditional forecasts is \emph{coherent$_1$} if, for
every finite
subset $\{i_1,\ldots,i_n\}\subseteq I$ and all real
$\alpha_1,\ldots,\alpha_n$, there exists no $\varepsilon>0$ such that
\[
\sum_{j=1}^n\alpha_jH_{i_j}(
\omega)\bigl[X_{i_j}(\omega)-P(X_{i_j}|H_{i_j})\bigr]
\geq\varepsilon
\]
for all $\omega\in\Omega$.

A collection of conditional forecasts is \emph{coherent$_2$} if
no sum of finitely many (Brier score) penalties
can be uniformly strictly dominated in the partition of states by the
sum of penalties from a rival set of forecasts for the same random
variables. That is, for every finite subset
$\{i_1,\ldots,i_n\}\subseteq I$, all alternative forecasts
$q_{i_1},\ldots,q_{i_n}$, and all positive $\alpha_1,\ldots,\alpha_n$,
there is no $\varepsilon>0$ such that\looseness=-1
\[
\sum_{j=1}^n\alpha_jH_{i_j}(
\omega)\bigl[X_{i_j}(\omega)-P(X_{i_j}|H_{i_j})
\bigr]^2 \geq\sum_{j=1}^n
\alpha_jH_{i_j}(\omega)\bigl[X_{i_j}(\omega
)-q_{i_j}\bigr]^2+\varepsilon
\]\looseness=0
for all $\omega$.
\end{definition}

De Finetti [(\citeyear{deFinetti1974}), pages~88--89] proved that a decision maker who wishes
to be both coherent$_1$ and coherent$_2$ must choose the same
forecasts for both purposes.
He also proved that the decision maker's coherent$_1$ forecasts are
represented by a finitely additive personal
probability, $P(\cdot)$, in the sense of Definition~\ref{def:fap} below.

If $P(H)=0$, then coherence$_1$ and coherence$_2$ place no
restrictions on $P(X|H)$ for bounded $X$.
Nevertheless, it is possible and useful to make certain intuitive
assumptions about conditional forecasts given events with 0
probability. In particular, Theorems \ref{thm:3} and \ref{thm:2008b}
of this paper assume that $P(\cdot|H)$ is a
finitely additive expectation (in the sense of
Definition~\ref{def:daniell} in Appendix \ref{app:lf}) satisfying
$P(X|H)=P(HX|H)$ for all $H$ and $X$. This assumption holds whenever
$P(H)>0$, and it captures the idea that $P(\cdot|H)$ is concentrated
on $H$.
De Finetti [(\citeyear{deFinetti1975}), Appendix~16] introduces an axiom that places a
similar requirement on conditional previsions. See Levi
(\citeyear{Levi1980}), Section~5.6, and \citet{regazzini1987} for other ways
to augment the coherence criteria of Definition~\ref{def:condcoh}
in
order to satisfy these
added requirements on conditional previsions given a null event.
Rather than adding such requirements to the definition of coherence, we\vadjust{\goodbreak}
prefer that individual agents who wish to adopt them do so as
explicit additional assumptions. Example~2 in the
supplemental article [Schervish, Seidenfeld and
Kadane (\citeyear{ssk2013})] illustrates our reason
for such a preference. In this way, our definition of coherence
is slightly weaker than that of de Finetti.

As an aside, the meaning of conditional expected value in the
finitely-additive theory differs from its meaning in the
countably-additive theory in this one major regard: In the
finitely-additive theory a conditional expectation can be specified given
an arbitrary nonempty event, regardless of whether that event has positive
probability. A conditional expectation of a bounded random variable
given an event with zero
probability is not defined uniquely in terms of unconditional
expectations, but \citet{dubins1975} shows that, in the finitely additive
theory, conditional expectations can be defined on the set of bounded
random variables so that they are
finitely additive expectations. In the countably-additive theory,
conditional expectation is defined twice: given events with positive
probability and given $\sigma$-fields. The two definitions match in a
well-defined way, and both provide uniquely defined conditional
expectations in
terms of unconditional expectations.

\begin{definition}\label{def:fap}
A probability $P(\cdot)$ is \emph{finitely additive} provided that,
when events $F$ and $G$ are disjoint, that is, when $F \cap G =
\varnothing$,
then $P(F \cup G) = P(F) + P(G)$.
A probability is \emph{countably additive} provided that when $F_{i}\ (i = 1, \ldots)$ is a denumerable sequence of pairwise disjoint
events, that is, when $F_{i} \cap F_{j} = \varnothing$ if $i \neq j$, then
$P(\bigcup_{i=1}^\infty F_{i}) = \sum_{i=1}^\infty P(F_{i})$.
We call a probability $P$ \emph{merely finitely additive} when $P$
is finitely but not countably additive. Likewise, then its
$P$-expectations are merely finitely additive.
\end{definition}

For each pair $X$ and $Y$ of random variables with finite previsions
(expectations), $P(X +Y) = P(X) + P(Y)$.
For countably additive expectations and countably many random
variables $\{X_{i}\}_{i=1}^\infty$, conditions under which
$P (\sum_{i=1}^\infty X_i ) = \sum_{i=1}^\infty P(X_i)$ can
be derived from various theorems such as the monotone convergence theorem,
the dominated convergence theorem, Fubini's theorem and Tonelli's theorem.

De Finetti (\citeyear{deFinetti1981}) recognized that coherence$_2$ (but not
coherence$_1$) provided an incentive compatible solution to the
problem of mechanism design for eliciting a coherent set of personal
probabilities. Specifically, Brier score is a \textit{strictly
proper scoring rule}, as defined here.

\begin{definition}\label{def:4} A \emph{scoring rule} for coherent
forecasts of a random
variable $X$ is a real-valued loss function $g$ with two real arguments:
a value of the random variable and a forecast $q$.
Let $\P_g$ be the collection of probability distributions such that
$P(X)$ is finite and $P[g(X,q)]$ is finite for at least one $q$.
We say that $g$ is \emph{proper} if, for every probability $P\in\P_g$,
$P[g(X,q)]$ is minimized (as a function of $q$) by $q=P(X)$.
If, in addition, only the quantity $q=P(X)$ minimizes expected score,
then the scoring rule is \emph{strictly proper}.\vadjust{\goodbreak}
\end{definition}

The following trivial result connects proper scoring rules with
conditional distributions.
%
\begin{proposition}\label{pro:src}
If $H$ is a nonempty event and $P(\cdot|H)$ is a probability distribution
then $P[g(X,q)|H]$ is (uniquely) minimized by $q=P(X|H)$ if $g$ is
(strictly) proper.
\end{proposition}

Some authors reserve the qualification strictly proper for scoring
rules that are designed to elicit an entire distribution, rather than
just the mean of a distribution.
[See Gneiting (\citeyear{gneiting11a}), who calls the
latter kind \emph{strictly consistent}.] For the remainder of this
paper, we
follow the language of Definition~\ref{def:4}, which matches the usage in
Gneiting (\citeyear{gneiting11b}).

We present some background on strictly proper scoring rules in
Section~\ref{sec:scoring}. Section~\ref{sec:extensions} gives our
main results.
We discuss propriety of scoring rules for infinitely many forecasts in
Section~\ref{sec:proper}.

\section{Background on strictly proper scoring rules}\label{sec:scoring}
In this section, we introduce a large class of strictly proper scoring
rules that we use as generalizations of Brier score. Associated with
this class, we introduce a third coherence concept that generalizes
coherence$_2$.

\begin{definition}
Let ${\cal C}$ be a class of strictly proper scoring rules. Let
$\{(X_i, H_i)\dvtx i\in I\}$ be a collection of random variable/nonempty
event pairs with corresponding conditional forecasts $\{p_i\dvtx i\in I\}$.
The forecasts are\break  \emph{coherent$_3$ relative to ${\cal C}$} if, for
every finite subset $\{i_{j}\dvtx j=1,\ldots,n\}\subseteq I$, every set of scoring rules
$\{g_j\}_{j=1}^n\subseteq{\cal C}$, and every set $\{q_j\}_{j=1}^n$ of
alternative forecasts, there is no $\varepsilon>0$ such that
\[
\sum_{j=1}^nH_{i_j}(
\omega)g_j\bigl(X_{i_j}(\omega),p_{i_j}\bigr)\geq
\sum_{j=1}^n H_{i_j}(
\omega)g_j\bigl(X_{i_j}(\omega),q_j\bigr)+
\varepsilon
\]
for all $\omega$.
That is, no sum of finitely many scores can be uniformly strictly
dominated by the sum of scores from rival forecasts.
\end{definition}
Coherence$_2$ is the special case of coherence$_3$ in which ${\cal C}$
consists solely of Brier score. The supplemental article
[Schervish, Seidenfeld and
Kadane (\citeyear{ssk2013})] includes a proof
that, if ${\cal C}$ consists of strictly proper
scoring rules of the form (\ref{eq:gdef}) below, then coherence$_3$
relative to ${\cal C}$ is equivalent to coherence$_1$.

The general form of scoring rule that we will consider is
%
\begin{equation}
\label{eq:gdef} g(x,q)=\cases{ %
\displaystyle\int
_x^q(v-x)\,d\lambda(v),&\quad $\mbox{if $x\leq q$,}$
\vspace*{2pt}\cr
\displaystyle\int_q^x(x-v)\,d\lambda(v),&\quad $\mbox{if $x>q$,}$}
\end{equation}
where $\lambda$ is a measure that is mutually absolutely continuous with
Lebesgue measure and is finite on every bounded interval.
It is helpful to rewrite (\ref{eq:gdef}) as
%
\begin{equation}
\label{eq:gsimp} g(x,q)=\int_q^x(x-v)\,d\lambda(v),
\end{equation}
using the convention that an integral whose limits are in the wrong
order equals the negative of the integral with the limits in the
correct order. Another interesting way to rewrite (\ref{eq:gdef}),
using the same convention, is
%
\begin{equation}
\label{eq:glin} g(x,q)=\lambda\bigl((q,x)\bigr)\bigl[x-r(x,q,\lambda)\bigr],
\end{equation}
where, for all $a$ and $b$,
%
\begin{equation}
\label{eq:ra2} r(a,b,\lambda)=\frac{\int_a^bv\,d\lambda(v)}{\lambda((a,b))}.
\end{equation}
An immediate consequence of (\ref{eq:glin}) is that, if $p$ and $q$
are real numbers, then
%
\begin{equation}
\label{eq:diff} g(x,q)-g(x,p)=\lambda\bigl((q,p)\bigr)\bigl[x-r(q,p,\lambda)
\bigr].
\end{equation}

The form (\ref{eq:gdef}) is suggested by equation (4.3) of \citet{savage71}.
Each such scoring rule is finite, nonnegative and continuous as a
function of $(x,q)$. If we wanted to
consider only countably additive distributions, we could use a larger
class of scoring rules by allowing $\lambda$ to be an infinite measure
supported on a bounded interval $(c_1,c_2)$. But this relaxation would
allow functions $g$ that are not strictly proper for natural classes
of finitely additive distributions. Example~\ref{exa:log} below
illustrates this point. Lemma~\ref{lem:sp} justifies the
use of (\ref{eq:gdef}) as the form of our scoring rules. The proofs of
all results in the body of the paper are given in Appendix \ref{app:aux}.

\begin{lemma}\label{lem:sp}
Let $g$ be a scoring rule of the form (\ref{eq:gdef}). Then $g$ is strictly
proper.
\end{lemma}
It follows from (\ref{eq:glin}) that, if $\lambda$ is a probability
measure with finite mean, then $\P_g$ from Definition~\ref{def:4} is
the class of all finitely additive distributions with finite mean because
$\lambda((q,x))$ and $\lambda((q,x))r(q,x,\lambda)$ are both bounded
functions of $q$ and $x$. Even if $\lambda$ is not a finite measure,
(\ref{eq:diff}) implies that, if $P[g(X,q_0)]$ is finite, then
$h(x,p)=g(x,p)-g(x,q_0)$ is linear in $x$ so that $h(x,p)$ is also
strictly proper with $\P_h$ equal to the class of all probabilities
with finite mean. For example, if $g(x,p)=(x-p)^2$, namely Brier score, then
$\P_g$ is the set of distributions with finite second moment. However,
$h(x,p)=(x-p)^2-x^2$ has $\P_h$ equal to the class of all probabilities
with finite mean.

Let $f(\cdot)$ denote the Radon--Nikodym derivative of
$\lambda$ with respect to Lebes\-gue measure.
Some familiar examples of strictly proper scoring rules are recovered
by setting $f$ equal to specific functions. Brier score corresponds
to $f(v)\equiv2$. Logarithmic score on the interval $(c_1,c_2)$ corresponds
to $f(v)=(c_2-c_1)/[(c_2-v)(v-c_1)]$, but the corresponding measure is
infinite on $(c_1,c_2)$. Hence, logarithmic score is not of the form
(\ref{eq:gdef}).
In addition, if $g$ is this logarithmic score, then $\P_g$ does not
include all finitely additive distributions that take values in the
bounded interval $(c_1,c_2)$, as the following example illustrates.
%
\begin{example}\label{exa:log}
Let $X$ be a random variable whose entire distribution
is agglutinated at $c_1$ from above. That is, let $P(X>c_1)=1$ and
$P(X<c_1+\varepsilon)=1$ for all $\varepsilon>0$. Let $g$ be the
logarithmic scoring rule that uses $f(v)$ from above.
Then $P(X)=c_1$, but $g(X(\omega),c_1)=\infty$ for all $\omega$,
which could not have finite mean even if we tried to extend the definition
of random variables to allow them to assume infinite values. On the
other hand, for $c_1<q<c_2$, the mean of $g(X,q)$ is
$\log[(c_2-c_1)/(c_2-q)]>0$, which decreases to 0 as $q$ decreases to
$c_1$, and is always finite. So, $\P_g$ is nonempty but does not
contain~$P$.
\end{example}

Some of our results rely on one or another condition that prevents the
$\lambda$ measures that determine the scoring rules from either being
too heavily concentrated on small sets or from being too different
from each other.
%
\begin{definition}
Let ${\cal C}=\{g_i\dvtx i\in I\}$ be a collection of strictly proper
scoring rules of the form (\ref{eq:gdef}) with corresponding measures
$\{\lambda_i\dvtx i\in I\}$.
\begin{longlist}[(ii)]
\item[(i)] Suppose that, for every $\varepsilon\geq0$, there exists
$\delta_\varepsilon>0$ such that for all $i\in I$ and all real $a<b$,
$\lambda_i((a,b))>\varepsilon$ implies $a+\delta_\varepsilon\leq
r(a,b,\lambda
_i)\leq
b-\delta_\varepsilon$.
Then we say that the collection ${\cal C}$ satisfies the \emph{uniform
spread condition}.
\item[(ii)] Suppose that, for every $\varepsilon>0$ and every $i\in I$,
there exists $\gamma_{i,\varepsilon}>0$ such that for all $j\in I$ and
all real $a<b$, $\lambda_i((a,b))\geq\varepsilon$ implies
$\lambda_j((a,b))\geq\gamma_{i,\varepsilon}$. Then
we say that the collection ${\cal C}$ satisfies the \emph{uniform
similarity condition}.
\end{longlist}
\end{definition}
The $r(a,b,\lambda)$ in (\ref{eq:ra2}) can be thought of as the mean
of the probability measure on the interval $(a,b)$ obtained by
normalizing $\lambda$ on the interval. The uniform spread condition
insures that, the $\lambda$ measures are spread out enough to keep the means
of the normalized measures on intervals far enough away from both endpoints.

The next result gives sufficient conditions for both the uniform similarity
and uniform spread conditions. It is easy to see that the conditions
are logically independent of each other.
%
\begin{lemma}\label{lem:ubf}
Let ${\cal C}=\{g_i\dvtx i\in I\}$ be a collection of strictly proper
scoring rules of the form (\ref{eq:gdef}) with corresponding measures
$\{\lambda_i\dvtx i\in I\}$ and corresponding Radon--Nikodym derivatives
$\{f_i\dvtx i\in i\}$ with respect to Lebesgue measure.
\begin{longlist}[(ii)]
\item[(i)] Assume that there exists $U<\infty$ such
that $f_i(v)\leq U$, for all $v$ and all $i\in I$. Then ${\cal C}$
satisfies the uniform spread condition.
\item[(ii)] Assume that for every $i\in I$, there exists $L_i>0$ such
that $f_j(v)/f_i(v)\geq L_i$, for all $v$ and all $j\in I$. Then ${\cal C}$
satisfies the uniform similarity condition.
\end{longlist}
\end{lemma}

As an example, suppose that each
$\lambda_i$ is $\alpha_i>0$ times Lebesgue measure. If the
$\alpha_i$ are bounded above, then ${\cal C}$
satisfies the uniform spread condition. If the $\alpha_i$ are bounded
away from 0, then ${\cal C}$ satisfies the uniform similarity
condition. These sets of measures correspond to
multiples of Brier score. There are collections that satisfy the uniform
spread condition without satisfying the conditions of part (i) of
Lemma~\ref{lem:ubf}. For example, let $f(v)=|v|^{-1/2}/2$ which is
not bounded above. For this~$f$, we have
$\lambda((a,b))=|\sqrt{|b|}-\sqrt{|a|}|$ if $0\notin(a,b)$, and
$\lambda((a,b))=\sqrt{|a|}+\sqrt{|b|}$ if $0\in(a,b)$.
So $\lambda((a,b))^2$ is no larger than two times the distance between
$a$ and~$b$. Also, $r(a,b,\lambda)$ is always at least $1/3$ of the way
from both $a$ and $b$. We can add the corresponding scoring rule to
any class that already satisfies the uniform spread condition by (if
necessary) lowering $\delta_\varepsilon$ to $\varepsilon^2/6$.

\section{Extensions to countably many options}\label{sec:extensions}
In Section~\ref{sub:dominance3}, we investigate when each sense of
coherence can be extended to allow combining countably many forecasts into
a single act by summing together their individual outcomes. In
Section~\ref{sub:dominance4}, we introduce the concept of \emph{conditional} forecasts and present results about the
combination of countably many coherent conditional forecasts.

\subsection{Dominance for countably many forecasts}\label{sub:dominance3}
Let $\{X_i\}_{i=1}^\infty$ be a countable set of random variables
with corresponding coherent$_1$ unconditional previsions $\{p_i\}
_{i=1}^\infty$.
Let $\{\alpha_i\}_{i=1}^\infty$ be a sequence of real numbers.
The decision maker's net loss
in state $\omega$, from adding the individual losses from the
fair options $\alpha_i[X_i(\omega)-p_i]$ is
%
\begin{equation}
\label{eq:infsumprev} \sum_{i=1}^\infty
\alpha_{i}\bigl[X_{i}(\omega) - p_{i}\bigr].
\end{equation}
Similarly, if the agent's prevision $p_i$ for $X_i$ is scored by the
strictly proper scoring rule $g_i$ for each $i$, the total score in
each state $\omega$ equals
\[
\sum_{i=1}^\infty g_i
\bigl(X_i(\omega), p_{i}\bigr).
\]
We assume that each of the two series above are convergent for all
$\omega\in\Omega$.

\begin{example}[(Combining countably many
forecasts)]\label{examp:1} De
Finetti [(\citeyear{deFinetti1949}), page~91] noted that when the decision
maker's personal probability is merely finitely additive, she/he
cannot always accept as fair the countable sum~(\ref{eq:infsumprev})
determined by coherent$_{1}$ forecasts. That sum may be uniformly
dominated by abstaining. Let $\Omega=
\{\omega_i\}_{i=1}^\infty$ be a countable state space. Let $W_i$ be
the indicator function for state $\omega_i\dvtx W_i(\omega) = 1$ if
$\omega
= \omega_i$ and $W_i(\omega) = 0$ if $\omega\neq\omega_i$.
Consider a collection of merely finitely additive coherent$_1$ forecasts
$P(W_i) = p_i \geq0$ where $\sum_{i=1}^\infty p_{i} = c < 1$. So
$P(\cdot)$ is not countably additive. With
$\alpha_i=1$, for all $i$, the loss from combining these
infinitely many forecasts into a single option is uniformly positive,
\[
\sum_{i=1}^\infty\alpha_i
\bigl[W_i(\omega) - p_i\bigr] =(1-c)>0.
\]
Hence, the
decision maker's alternative to abstain, with constant loss $0$,
uniformly strictly dominates this infinite combination of fair options.

If, on the other hand, the decision maker's personal probability $P$
is countably additive, then $c=1$. For arbitrary
$\{\alpha_i\}_{i=1}^\infty$ such that $d=\sum_{i=1}^\infty\alpha_ip_i$
is defined and finite, the sum of losses is
%
\begin{equation}
\label{eq:gca} \sum_{i=1}^\infty
\alpha_i\bigl[W_i(\omega) - p_i\bigr] =
\alpha_{i(\omega)}-d,
\end{equation}
where $i(\omega)$ is the unique $i$ such that $W_i(\omega)=1$.
Because $c=1$, there is at least one $\alpha_i\leq d$ and at least one
$\alpha_i\geq d$, hence (\ref{eq:gca}) must be nonpositive for at
least one~$i$, and abstaining does not uniformly strictly dominate.
\end{example}

Next, we focus on the parallel question whether a coherent$_3$ set
of forecasts remains undominated when strictly proper scores for
countably many forecasts are summed together. Some conditions will
be needed in order to avoid $\infty-\infty$ arising in the
calculations, and these are stated precisely in the theorems.
The principal difference between dominance for infinite sums of
forecasts and dominance for infinite sums of strictly proper scores is
expressed by the following result.

\begin{theorem}\label{thm:1}
Let ${\cal C}$ be a collection of strictly proper scoring rules of
the form~(\ref{eq:gdef}) that satisfies the uniform spread condition.
Let $P$ be a coherent$_3$ prevision defined over a collection $\D$ of
random variables that contains all of the random variables mentioned
in the statement of this theorem. Let $\{X_i\}_{i=1}^\infty$ be
random variables in $\D$ with
coherent$_3$ forecasts $P(X_i) = p_{i}$ for $i=1,2,\ldots.$ Assume
that the forecast for $X_i$ will be scored by a scoring rule
$g_i\in{\cal C}$ for each $i$. Finally, assume that
%
\begin{eqnarray}
\label{eq9} P \Biggl[\sum_{i=1}^\infty\mid
X_{i}-p_i \mid \Biggr] &=& V < \infty\quad\mbox{and}
\\
\label{eq2} P \Biggl[\sum_{i=1}^\infty
g_i(X_i,p_i) \Biggr] &= & W < \infty.
\end{eqnarray}
There does not exist a rival set of forecasts $\{q_i\}_{i=1}^\infty$
such that, for all $\omega\in\Omega$,
%
\begin{equation}
\label{eq:dom} \sum_{i=1}^\infty
g_i\bigl(X_i(\omega),p_i\bigr)>\sum
_{i=1}^\infty g_i
\bigl(X_i(\omega),q_i\bigr).
\end{equation}
\end{theorem}

Theorem~\ref{thm:1} asserts conditions under which infinite sums of
strictly proper scores, with coherent$_3$ forecasts $\{p_{i}\}
_{i=1}^\infty$ for
$\{X_{i}\}_{i=1}^\infty$, have no rival forecasts that
simply dominate, let alone uniformly strictly
dominate $\{p_i\}_{i=1}^\infty$. That is, even countably many
unconditional coherent$_3$ forecasts cannot be simply dominated
under the conditions of Theorem~\ref{thm:1}.

\begin{example}[(Example~\ref{examp:1} continued)]\label{examp:3}
Recall that
$\Omega= \{\omega_{i}\dvtx i = 1, \ldots\}$ is a countable space.
Consider the special case in which $P$ is a purely finitely
additive probability satisfying
$P(\{\omega_{i}\}) = p_{i} = 0$, for all $i$. So, $\sum_{i=1}^\infty
p_{i} = 0 < 1$, and $c=0$ in the notation of
Example~\ref{examp:1}. As before, let $W_{i}$ ($i = 1, \ldots$) be
the indicator functions for the states in $\Omega$. So $P(W_i)
= p_{i} = 0$ and combining the losses $W_{i} - p_{i} =W_{i}$, for
$i = 1, \ldots$ results in a uniform
sure-loss of $1$. But this example, with each $g_i$ equal to Brier
score times $\alpha_i>0$, satisfies
the conditions of Theorem~\ref{thm:1}, if the $\alpha_i$ are
bounded above. That is, there are no
rival forecasts $\{q_i\}_{i=1}^\infty$ for the
$\{W_{i}\}_{i=1}^\infty$ that simply dominate the
forecasts $\{p_i\}_{i=1}^\infty$ by weighted sum of Brier scores, let alone
uniformly strictly dominating these forecasts.
We can illustrate the conclusion of Theorem~\ref{thm:1} directly in
this example. The weighted sum of Brier
scores for the $p_i$ forecasts is\looseness=-1
\[
S(\omega)=\sum_{i=1}^\infty
\alpha_iW_i(\omega)^2\leq\sup
_i\alpha_i.
\]\looseness=0
Let $\{q_i\}_{i=1}^\infty$ be a rival set of
forecasts with $q_i\ne p_i$ for at least one $i$. The corresponding
weighted sum of Brier scores is
%
\begin{equation}
\label{eq:exampsum} \sum_{i=1}^\infty
\alpha_i\bigl[W_i(\omega)-q_i
\bigr]^2=S(\omega) -2\sum_{i=1}^\infty
\alpha_iq_iW_i(\omega)+\sum
_{i=1}^\infty\alpha_iq_i^2.
\end{equation}
Let $d=\sum_{i=1}^\infty\alpha_iq_i^2$, which must be strictly greater
than 0.
Define $i(\omega)$ to be the unique value of $i$ such that $W_i(\omega)=1$.
The right-hand side of (\ref{eq:exampsum}) can then be written as
$S(\omega)-2\alpha_{i(\omega)}q_{i(\omega)}+d$. If $d=\infty$, then
the rival forecasts clearly fail to dominate the original forecasts.
If $d<\infty$, then $\lim_{i\rightarrow\infty}\sqrt{\alpha_i}q_i=0$.
Because the $\alpha_i$ themselves are bounded, it follows that all but finitely
many $\alpha_i|q_i|$ are less than $d/2$. For each $\omega$ such that
$\alpha_{i(\omega)}|q_{i(\omega)}|<d/2$, we have the weighted sum of
Brier scores displayed in (\ref{eq:exampsum})
strictly greater than $S(\omega)$, hence the rival forecasts do not
dominate the original forecasts.
\end{example}

Theorem~\ref{thm:1}, as illustrated by Example~\ref{examp:3}, shows
that the modified decision problem in de Finetti's prevision game---modified
to include infinite sums of betting outcomes---is not
isomorphic to the modified forecasting problem under strictly proper
scoring rules---modified to include infinite sums of scores. In particular,
abstaining from betting, which is the alternative that uniformly
dominates the losses for coherence$_1$, is not an available alternative
under forecasting with strictly proper scores. In summary, the two
criteria, coherence$_1$ and coherence$_3$ behave differently when
probability is merely finitely additive and we try to combine
countably many forecasts.\looseness=1

We conclude this section with an example to show why we assume that
the class of scoring rules satisfies the uniform spread condition in
Theorem~\ref{thm:1}.

\begin{example}
This example satisfies all of the conditions of Theorem~\ref{thm:1}
except that the class of scoring rules fails
the uniform spread condition. We show that the conclusion to
Theorem~\ref{thm:1} also fails. For each integer $i\geq1$, let
$\alpha_i=2^{-i-1}$, and define
\[
f_i(v)=\cases{ %
2,&\quad $\mbox{if $v\leq
\alpha_i$,}$
\vspace*{2pt}\cr
\displaystyle\frac{2}{\alpha_i}, &\quad
$\mbox{if $v>\alpha_i$.}$}
\]
Let $\lambda_i$ be the measure whose Radon--Nikodym derivative with
respect to Lebesgue measure is $f_i$, and define $g_i$ by
(\ref{eq:gdef}) using $\lambda=\lambda_i$. The form of $g_i$ is as follows:
\begin{eqnarray*}
g_i(x,q)=\cases{ %
(x-q)^2,&\quad
$\mbox{if $x, q\leq\alpha_i$,}$
\vspace*{2pt}\cr
\displaystyle(x-\alpha_i)^2+\frac{1}{\alpha_i} (q-
\alpha_i)^2+\frac{2}{\alpha
_i}(q-\alpha _i) (
\alpha_i-x), &\quad $\mbox{if $x\leq\alpha_i\leq q$,}$
\vspace*{2pt}\cr
\displaystyle\frac{1}{\alpha_i}(x-\alpha_i)^2+(q-
\alpha_i)^2+2(\alpha _i-q) (x-
\alpha_i) ,&\quad$\mbox{if $q\leq\alpha_i\leq x$,}$
\vspace*{2pt}\cr
\displaystyle\frac{1}{\alpha_i}(x-q)^2,&\quad $\mbox{if $\alpha_i\leq x, q$.}$}
\end{eqnarray*}
These scoring rules fail
the uniform spread condition because arbitrarily short intervals with
both endpoints positive have arbitrarily large $\lambda_i$ measure as
$i$ increases. Let $\{A_i\}_{i=1}^\infty$ be a partition of the real
line, and let $P(\cdot)$ be a finitely additive probability such that
$P(A_i)=0$ for all $i$. For each integer $i\geq1$, let $p_i=2^{-i}$ and
$q_i=2^{-i-1}$, and define
\[
X_i(\omega)=\cases{ %
p_i,&\quad
$\mbox{if $\omega\in A_i^C$,}$
\vspace*{2pt}\cr
q_i-1,&\quad $\mbox{if $\omega\in A_i$.}$}
\]
It follows that $P(X_i)=p_i$ for all $i$, and
\[
P \Biggl[\sum_{i=1}^\infty|X_i-p_i|
\Biggr]=P \Biggl[ \sum_{i=1}^\infty
A_i|q_i-1-p_i| \Biggr]=1,
\]
so that (\ref{eq9}) holds. Next, compute the various scores:
\begin{eqnarray*}
&&g_i\bigl(X_i(\omega),p_i
\bigr)\\
&&\qquad=A_i(\omega) \biggl\{ (q_i-1-
\alpha_i)^2+\frac{1}{\alpha_i}(p_i-
\alpha_i)^2+ \frac{2}{\alpha_i}(p_i-
\alpha_i) (\alpha_i-q_i+1) \biggr\}
\\
&&\qquad=A_i(\omega)\bigl[3+2^{-i-1}\bigr],
\\
&&g_i\bigl(X_i(\omega),q_i
\bigr)\\
&&\qquad=A_i^C(\omega)\frac{1}{\alpha_i}(p_i-q_i)^2
+A_i(\omega) (q_i-1-q_i)^2
\\
&&\qquad=A_i^C(\omega)2^{-i-1}+A_i(
\omega).
\end{eqnarray*}
Define $i(\omega)=i$ for that unique $i$ such that $\omega\in A_i$.
When we sum up the scores for the forecasts $\{p_i\}_{i=1}^\infty$, we get
\[
\sum_{i=1}^\infty g_i
\bigl(X_i(\omega),p_i\bigr)=\sum
_{i=1}^\infty A_i(\omega)
\bigl[3+2^{-i-1}\bigr]=3+2^{-i(\omega)-1}.
\]
It follows that $P(\sum_{i=1}^\infty g_i(X_i,p_i))=3$, so (\ref{eq2}) holds.
The sum of the $\{q_i\}_{i=1}^\infty$ scores is
\[
\sum_{i=1}^\infty g_i
\bigl(X_i(\omega),q_i\bigr)=\sum
_{i=1}^\infty \bigl[A_i^C(
\omega)2^{-i-1} +A_i(\omega)\bigr]=1.5-2^{-i(\omega)-1}.
\]

Finally, compute the difference in total scores:
\[
\sum_{i=1}^\infty g_i
\bigl(X_i(\omega),p_i\bigr)-\sum
_{i=1}^\infty g_i\bigl(X_i(
\omega),q_i\bigr) =1.5+2^{-i(\omega)}>1.5,
\]
hence the scores of the $\{p_i\}_{i=1}^\infty$ forecasts are uniformly
strictly dominated by the scores of a set of rival forecasts.
\end{example}

\subsection{Dominance for countable sums of conditional
forecasts}\label{sub:dominance4}
Definition~\ref{def:condcoh} allows mixing conditional forecasts with
unconditional forecasts by setting $P(X|H)=P(X)$ whenever $H=\Omega$.
De Finetti showed that, if $P(X)$, $P(X|H)$ and $P(HX)$ are all
specified, a necessary condition for coherence$_1$ is
that
%
\begin{equation}
\label{eq:cond} P(HX)=P(H)P(X|H),
\end{equation}
so that $P(X|H)$ is the usual conditional expected value of $X$ given
$H$ whenever $P(H)>0$. For this reason, conditional forecasts are
often called \emph{conditional expectations}.

The concept of \textit{conglomerability} plays a central role in our
results about coherence for combining countably many conditional
forecasts. Conglomerability in a partition $\pi=
\{H_j\dvtx j\in J \}$ of conditional expectations $P(\cdot
\mid H_j)$ over a class $\D$ of random variables $X$ is the
requirement that the unconditional expectation of each $X\in\D$
lies within the range of its conditional expectations given
elements of $\pi$.

\begin{definition}\label{def:nonc} Let $P$ be a finitely additive
prevision on a set $\D$ of random variables,
and let $\pi=\{H_j\dvtx j\in J\}$ be a partition of $\Omega$ such that
conditional prevision $P(\cdot|H_j)$ has been defined for all
$j$. If, for each $X\in\D$,
\[
\inf_{j \in J} P(X \mid H_{j}) \leq P(X) \leq\sup
_{j \in J} P(X \mid H_{j}),
\]
then $P$ is
\emph{conglomerable in the partition $\pi$ with respect to $\D$}.
Otherwise, $P$ is \emph{nonconglomerable} in $\pi$ with respect to
$\D$.
\end{definition}


If a decision maker's coherent$_1$ or coherent$_3$
forecasts fail conglomerability in a partition $\pi$,
Theorem~\ref{thm:2} below shows
there exist countably many conditional forecasts that are uniformly
strictly dominated.

On the other hand, if the decision maker's previsions for random
variables satisfy a condition (see Definition~\ref{def:ltp}) similar
to being conglomerable in $\pi$, Theorem~\ref{thm:3} below
establishes that no countable set of forecasts,
conditional on elements of $\pi$, can be uniformly strictly dominated.
What we mean by ``similar'' is explained in Section~\ref{sec:cdltp} below.

\begin{theorem}\label{thm:2}
Let $P$ be a finitely additive prevision, and let $\D$ be a set of
random variables. Let $\pi= \{H_j\}_{j=1}^\infty$ be a denumerable
partition and let $P(\cdot\mid H_j)$ be the corresponding conditional
previsions associated with $P$. Let ${\cal C}$ be a collection of strictly
proper scoring rules of the form (\ref{eq:gdef}) that satisfies the
uniform similarity condition.
Assume that the conditional previsions $P(\cdot\mid H_j)$
are nonconglomerable in $\pi$ with respect to $\D$. Then there exists
a random variable $X\in\D$ with $p_X=P(X)$ and $p_j=P(X|H_j)$ for all
$j$ such that

{(2.1)} the countable sum
\[
\alpha_0(X-p_X)+\sum_{j=1}^\infty
\alpha_jH_j(X-p_j),
\]
of individually fair options is uniformly strictly dominated by
abstaining, and

{(2.2)}
if the forecast for $X$ is scored by $g_0\in{\cal C}$ and the
conditional forecast for $X$ given $H_j$ is scored by $g_j\in{\cal C}$
for $j=1,2,\ldots,$ then the sum of the scores,
\[
g_0\bigl(X(\omega),p_X\bigr)+\sum
_{j=1}^\infty H_j(\omega)g_j
\bigl(X(\omega),p_j\bigr),
\]
is uniformly strictly dominated by the sum of scores from a rival set
of forecasts.
\end{theorem}

We illustrate Theorem~\ref{thm:2} with an example of
nonconglomerability due to \mbox{\citet{dubins1975}}. This example is
illuminating as the conditional probabilities do not involve
conditioning on null events.

\begin{example}\label{examp:2}
Let $\Omega= \{\omega_{ij}\dvtx i = 1,2; j = 1, \ldots\}$. Let $F =
\{\omega_{2j}, j = 1, \ldots\}$ and let $H_{j} = \{\omega_{1j},
\omega_{2j}\}$. Define a merely finitely additive probability $P$ so
that $P(\{\omega_{1j}\}) = 0, P(\{\omega_{2j}\}) = 2^{-(j+1)}$ for $j
= 1, \ldots,$ and let $P(F) = p_{F} = 1/2$. Note that $P(H_{j}) =
2^{-(j+1)} > 0$, so $P(F \mid
H_{j}) = 1 = p_{j}$ is well defined by the multiplication rule for
conditional probability. Evidently, the conditional probabilities
$\{P(F \mid H_{j})\}_{j=1}^\infty$ are nonconglomerable in $\pi$
since $P(F) = 1/2$ whereas $P(F \mid H_{j}) = 1$ for all $j$.

For {(2.1)}, Consider the fair options $\alpha_{j}H_{j}(F -
p_{j})$ for $j = 1, \ldots$ and $\alpha_F(F - p_F)$.
Choose $\alpha_j = 1$ and $\alpha_F=-1$. Then
\begin{eqnarray*}
&&\Biggl[-\bigl(F(\omega) - p_F\bigr) + \sum
_{j=1}^\infty H_j(\omega)\bigl[F(\omega) -
p_j\bigr] \Biggr] \\
&&\qquad= \cases{ %
0.5 - 1.0 = -0.5,&\quad $\mbox{if $\omega\notin F$,}$
\vspace*{2pt}\cr
-0.5 + 0.0 = -0.5,&\quad $\mbox{if $\omega\in F$.}$}
\end{eqnarray*}

Hence, these infinitely many individually fair options
are not collectively fair when taken together. Their sum is
uniformly strictly dominated by $0$ in $\Omega$, corresponding to the
option to abstain from betting.

Regarding {(2.2)}, unlike the situation with
Theorem~\ref{thm:1} involving countably many unconditional forecasts,
the sum of Brier scores from these conditional forecasts are
uniformly strictly dominated. In particular, the sum of Brier scores
for these forecasts is
\begin{eqnarray*}
&&\bigl(F(\omega) - p_F\bigr)^2 + \sum
_{j=1}^\infty H_j(\omega)\bigl[F(
\omega)-p_j\bigr]^2\\
&&\qquad= \cases{
0.25 + 1.00 = 1.25,&\quad $\mbox{if $\omega\notin F$,}$
\vspace*{2pt}\cr
0.25 + 0.00 = 0.25,&\quad $\mbox{if $\omega\in F$.}$}
\end{eqnarray*}
Consider the rival forecasts $Q(F \mid H_{j}) = 0.75 = q_{j}$ and
$Q(F) = 0.75 = q_F$. These correspond to the
countably additive probability $Q(\{\omega_{1j}\}) = 0.25 \times
2^{-j}$ and $Q(\{\omega_{2j}\}) = 0.75 \times2^{-j}$ for $j = 1,
\ldots.$ Then the combined Brier score from these countably many
rival forecasts is
\begin{eqnarray*}
&&\bigl(F(\omega)-q_F\bigr)^2+\sum
_{j=1}^\infty H_j(\omega)\bigl[F(
\omega)-q_j\bigr]^2 \\
&&\qquad=\cases{ %
9/16 + 9/16 = 1.125,&\quad $\mbox{if $\omega\notin F$,}$
\vspace*{2pt}\cr
1/16 + 1/16 = 0.125,&\quad $\mbox{if $\omega\in F$,}$}
\end{eqnarray*}
which is 0.125 less than the sum of the Brier scores of the original forecasts.
\end{example}

We offer one more example to show why we assume that the class of
scoring rules satisfies the uniform similarity condition in
Theorem~\ref{thm:2}.\vspace*{-1pt}
%
\begin{example}[(Example~\ref{examp:2} continued)]
Recall that we have a partition $\pi=\{H_j\}_{j=1}^\infty$ and an
event $F$ with $p_F=P(F)=0.5$ and $p_j=P(F|H_j)=1$ for all~$j$.
Let the unconditional forecast for $F$ be scored by Brier score, and let
the conditional forecast for $F$ given $H_j$ be scored by
$2^{-j-1}$ times Brier score. These scoring rules fail the
uniform similarity condition. We show that the conclusion to
Theorem~\ref{thm:2} fails. Specifically, we show that there is no
rival set of forecasts $q_F$ for $F$ and $q_j$ for $H_j$
($j=1,2,\ldots$) whose sum of scores uniformly strictly dominates the
original forecasts.

The total of the scores for the original forecasts is
%
\begin{equation}
\label{eq:pscore} \frac{1}{4}+\sum_{j=1}^\infty
H_j(\omega)2^{-j-1}\bigl[1-F(\omega)\bigr]^2.
\end{equation}
Consider an arbitrary rival set of forecasts with $q_F$ for $F$ and
$q_j$ for $F$ conditional on $H_j$. The sum of the scores for the rival
forecasts is
%
\begin{equation}
\label{eq:qscore} \bigl[q_F-F(\omega)\bigr]^2+\sum
_{j=1}^\infty H_j(\omega)2^{-j-1}
\bigl[q_j-F(\omega)\bigr]^2.
\end{equation}
Let $i(\omega)=j$ when $\omega\in H_j$. Then the difference
(\ref{eq:pscore}) minus (\ref{eq:qscore}) is
%
\begin{equation}
\label{eq:scorediff} \tfrac{1}{4}-\bigl[q_F-F(\omega)
\bigr]^2 +2^{-i(\omega)-1} \bigl(\bigl[1-F(\omega)
\bigr]^2-\bigl[q_j-F(\omega)\bigr]^2 \bigr).
\end{equation}
If $q_F=0.5$, then (\ref{eq:scorediff}) becomes
%
\begin{equation}
\label{eq:scorediff2} 2^{-i(\omega)}(1-q_{i(\omega)}) \biggl[\frac{1+q_{i(\omega
)}}{2}-F(
\omega ) \biggr].
\end{equation}
If there exists $\omega$ such that $q_{i(\omega)}\geq1$, (\ref
{eq:scorediff2})
is nonpositive, and the rival forecasts do not strictly dominate. If
all $q_{i(\omega)}<1$, (\ref{eq:scorediff2}) is negative for all
$\omega\in F$, and there is no dominance. If $q_F\ne0.5$, then
(\ref{eq:scorediff}) is at most
%
\begin{equation}
\label{eq:scorediff3} \tfrac{1}{4}-\bigl[q_F-F(\omega)
\bigr]^2+2^{-i(\omega)}.
\end{equation}
No matter what $q_F\ne0.5$ we pick, either $(q_F-1)^2$ or $(q_F-0)^2$ is
greater than $1/4$. Let $\delta=1/4-\max\{[q_F-1]^2,q_F^2\}$. For
$j>-\log_2(\delta)$, (\ref{eq:scorediff3}) is negative either for all
$\omega\in F\cap H_j$ or all $\omega\in F^C\cap H_j$. So, there is no
dominance.\vspace*{-1pt}
\end{example}

Last, we establish conditions under which combining strictly proper
scores from countably many conditional forecasts given elements of a
partition, or combining
the losses from countably many fair options based on those forecasts,
does not result in a uniform sure loss. A definition is useful first.

\begin{definition}\label{def:ltp}
Let $P$ be a finitely additive prevision on a set $\D$ of
random variables,
and let $\pi=\{H_j\dvtx j\in J\}$ be a\vadjust{\goodbreak} partition of $\Omega$ such that
conditional prevision $P(\cdot|H_j)$ has been defined for all
$j$. For each random variable $X\in\D$, we let $P(X|\pi)$ denote
the random variable $Y$ defined by $Y(\omega)=P(X|H_j)$ for all
$\omega\in H_j$ and all $j$. We say that $P$ \emph{satisfies the law
of total previsions in $\pi$ with respect to $\D$} provided that
for each random variable $X\in\D$, $P(X) = P[P(X \mid\pi)]$.
\end{definition}

\begin{theorem}\label{thm:3}
Let $P$ be a finitely additive prevision, and let $\D$ be a set of
random variables such that $P$ satisfies the law of
total previsions in $\pi= \{H_j\}_{j=1}^\infty$ with respect to $\D$.
Let $X\in\D$ be a random variable with
finite prevision $p_X=P(X)$ and finite conditional prevision
$p_j=P(X|H_j)$ given each $H_j$. Assume that $P(\cdot|H_j)$ is a finitely
additive expectation (in the sense of Definition~\ref{def:daniell})
that satisfies $P(X|H_j)=P(H_jX|H_j)$ for every $j$. Let ${\cal C}$ be
a collection of
strictly proper scoring rules.

{(3.1)}
Let $\{\alpha_j\}_{j=0}^\infty$ be real numbers. The sum of losses
%
\begin{equation}
\label{eq18} \alpha_0\bigl(X(\omega) - p_X\bigr) +
\sum_{j=1}^\infty\alpha_jH_j(
\omega)\bigl[X(\omega) - p_j\bigr],
\end{equation}
is not uniformly strictly dominated by abstaining.

{(3.2)} Let $g_0,g_1,\ldots$ be elements of ${\cal C}$. There is
no rival set of forecasts that uniformly strictly dominates the sum of scores
%
\begin{equation}
\label{eq24} g_0\bigl(X(\omega),p_X\bigr)+\sum
_{j=1}^\infty H_j(
\omega)g_j\bigl(X(\omega),p_j\bigr).
\end{equation}
\end{theorem}
%
\subsection{Conglomerability, disintegrability and the law of total
previsions}\label{sec:cdltp}

We claimed earlier that the law of total previsions in a partition
$\pi$ is similar to conglomerability in $\pi$. The claim begins with
a result of \citet{dubins1975}.
Dubins defines \emph{conglomerability in partition} $\pi$ of a finitely
additive prevision $P$ by the requirement that, for all bounded random
variables $X$,
\[
\mbox{if } \forall H \in\pi P(X |H) \geq0, \mbox{then } P(X) \geq0.
\]
Dubins' definition of conglomerability in $\pi$ is equivalent to
Definition~\ref{def:nonc} with respect to the set of all bounded
random variables. However, for a set $\D$ that includes unbounded random
variables and/or does not include all bounded random variables,
the two definitions are not equivalent without further assumptions.
Definition~\ref{def:nonc} is based on the definition given by de Finetti
[(\citeyear{deFinetti1974}), Section~4.7], which generalizes to unbounded random variables
more easily.

\citet{dubins1975} also defines \emph{disintegrability of $P$ in partition}
$\pi$
by the requirement that, for every bounded random variable $X$,
\[
P(X) = \int P(X |h) \,dP(h),
\]
where the finitely additive integral is as developed by Dunford and
Schwartz [(\citeyear{dunford-schwartz1958}), Chapter III].
Moreover, he establishes that conglomerability and disintegrability
in $\pi$ are equivalent for the class of bounded random
variables.\looseness=1

The law of total previsions in Definition~\ref{def:ltp}, with respect
to the set of all bounded random variables, is equivalent to
disintegrability in Dubins' sense, but not necessarily for sets that either
include some unbounded random variables or fail to include some bounded
random variables. In addition, not all real-valued coherent$_1$
previsions admit an integral
representation in the sense of Dunford and Schwartz for sets that
include unbounded random variables. For discussion of
the problem and related issues, see \citet{berti-etal2001};
Berti and Rigo
(\citeyear{berti-rigo1992,berti-rigo2000,berti-rigo2002}); Schervish, Seidenfeld,
and Kadane (\citeyear{ssk2008a}) and
\citet{ssk2009a}. As described in
Appendix \ref{app:lf}, we use a definition of
finitely additive integral that is a natural extension of coherent$_1$
prevision. In this way, the law of total previsions
extends Dubins' definition of disintegrability from bounded to
unbounded random variables without introducing the technical details
of Dunford and Schwartz. Finally, Theorem~1 of \citet{ssk2008a} gives conditions under which
conglomerability (Definition~\ref{def:nonc}) is equivalent to the law
of total previsions. The following is a translation of that result into
the notation and terminology of the present paper.
%
\begin{theorem}\label{thm:2008b}
Let $P$ be a finitely additive prevision on a set $\D$ of
random variables. Let $\pi= \{H_j\}_{j=1}^\infty$ be a denumerable
partition and let $P(\cdot\mid H_j)$ be the corresponding conditional
previsions associated with $P$. Assume that, for all $j$,
$P(\cdot|H_j)$ is a finitely additive expectation on $\D$. Also assume
that, for all $X\in\D$:
\begin{itemize}
\item$P(X)$ is finite,
\item$P(X|H_j)$ is finite for all $j$,
\item$H_jX\in\D$ for all $j$,
\item$P(H_jX|H_j)=P(X|H_j)$ for all $j$, and
\item$X-Y\in\D$, where $Y$ is defined (in terms of $X$) in
Definition~\ref{def:ltp}.
\end{itemize}
Then $P$ is conglomerable in $\pi$ with
respect to $\D$ if and only if $P$ satisfies the law of total
previsions in $\pi$ with respect to $\D$.
\end{theorem}

Under the conditions of Theorem~\ref{thm:2008b},
Theorems \ref{thm:2} and \ref{thm:3} show that, when the conditioning
events form a countable partition $\pi$, coherence$_1$ and
coherence$_3$ behave the same when extended to include, respectively,
the countable sum of individually fair options, and
the total of strictly proper scores from the forecasts. If and
only if these
coherent quantities are based on conditional expectations that are
conglomerable in $\pi$, then no failures of the Dominance Principle
result by combining infinitely many of them.

\citet{SchervishSeidenfeldandKadane1984} show that each merely
finitely additive probability fails to be conglomerable in some
countable partition. But each countably additive probability
has expectations that are conglomerable in each countable partition.
Thus, the conjunction of Theorems \ref{thm:1},
\ref{thm:2} and \ref{thm:3} identifies where the debate whether personal
probability may be merely finitely additive runs up against the debate
whether to extend either coherence criterion in order to apply it with
countable combinations of quantities.
We arrive at the following conclusions:

\begin{itemize}
\item Unless unconditional coherent$_{1}$ forecasts arise from a
countably additive probability,
combining countably many unconditional coherent$_{1}$ forecasts
into a single option may be dominated by abstaining.

\item However, under the conditions of Theorem~\ref{thm:1}, strictly
proper scoring rules are not similarly affected. The scores from
countably many coherent$_3$ unconditional
forecasts may be summed together without leading to a violation of
the Dominance Principle.

\item Unless conditional forecasts arise from
a set of conglomerable conditional probabilities, the
Dominance Principle does not allow combining countably many of these
quantities into a single option. Hence, only countably additive
conditional probabilities satisfy the Dominance Principle when an
arbitrary countable set of conditional quantities are summed
together.
\end{itemize}

\section{Incentive compatible elicitation of infinitely many forecasts
using strict\-ly proper scoring rules}\label{sec:proper}
Scoring an agent based on the values of the fair gambles constructed
from coherent$_1$ forecasts, is not proper. Because of the
presence of the opponent in the game, who gets to choose whether to
buy or to sell the random variable $X$ at the decision maker's
announced price, the decision maker faces a strategic choice of
pricing. For example, if the decision maker suspects that the
opponent's fair price, $Q(X)$, is greater than his own, $P(X)$, then
it pays to inflate the announced price and to offer the opponent,
for example, $R(X) = [P(X) + Q(X)]/2$, rather than offering $P(X)$.
Thus, the forecast-game as de Finetti defined it for coherence$_1$
is not incentive compatible for eliciting the decision maker's fair
prices.

With a finite set of forecasts and a strictly proper scoring rule
for each one, using the finite sum of the scores as the score
for the finite set preserves strict propriety. That is, with the sum
of strictly proper scores as the score for the finite set, a coherent forecaster
minimizes the expected sum of scores by minimizing each one, and this
solution is unique.

Here, we report what happens to the propriety of strictly proper scores
in each
of the three settings of the three theorems presented in
Section~\ref{sec:extensions}. That is, we answer the question whether
or not, in each of these three settings, the coherent forecaster
minimizes expected score for the infinite sum of strictly proper scores by
announcing her/his coherent forecast for each of the infinitely many
variables. These findings are corollaries to the respective
theorems.

\begin{corollary}
Under the assumptions of Theorem~\ref{thm:1}, the infinite sum of
scores applied to the infinite set of forecasts $\{p_i\}_{i=1}^\infty$
is a
strictly proper scoring rule.
\end{corollary}

\begin{corollary}
Under the assumptions used for {(2.2)} of Theorem~\ref{thm:2},
namely when
the conditional probabilities $P(F \mid H_j) = p_j$ are
nonconglomerable in $\pi$, then the infinite
sum of strictly proper scores applied to the infinite set of conditional
forecasts $\{p_j\}_{j=1}^\infty$ is not proper.
\end{corollary}

\begin{corollary}
Under the assumptions used to establish {(3.2)} of Theorem~\ref{thm:3},
namely that $P$ satisfies the law of total previsions in $\pi$, the
infinite sum of
strictly proper scores applied to the infinite set of conditional forecasts
$\{p_j\}_{j=1}^\infty$ is a proper scoring rule.
\end{corollary}

Thus, these results about the propriety of infinite sums of strictly proper
scores parallel the respective results about extending coherence$_3$
to allow infinite sums of scores.

\section{Summary}\label{sec:summary}
We study how two different coherence criteria
behave with respect to a Dominance Principle when countable
collections of random variables are included. Theorem~\ref{thm:1}
shows that, in
contrast with fair prices for coherence$_1$, when strictly proper
scores from infinitely many unconditional forecasts are summed together
there are no new failures of the Dominance Principle for
coherence$_3$. That is, if an infinite set of probabilistic forecasts
$\{p_{i}\}_{i=1}^\infty$ are even simply dominated by some rival
forecast scheme
$\{q_{i}\}_{i=1}^\infty$ in total score, then the $\{p_{i}\}
_{i=1}^\infty$ are not
coherent$_3$, that is, some finite subset of them is uniformly strictly
dominated in total score. However, because each merely finitely
additive probability fails to be conglomerable in some denumerable
partition, in the light of Theorem~\ref{thm:2}, neither of the two
coherence criteria discussed here may be relaxed in order to apply the
Dominance Principle with infinite combinations of conditional options.
Merely finitely additive probabilities then would become
\textit{incoherent}.

Specifically, the conjunction of
Theorems \ref{thm:1}--\ref{thm:2008b} shows that it matters only in
cases that involve
nonconglomerability whether incoherence$_3$ is established using
scores from a finite rather than from an infinite
combination of forecasts. In that one respect, we think
coherence$_3$ constitutes an improved version of the concept of
coherence. Coherence$_1$ applied to a merely finitely
additive probability leads to failures of the Dominance Principle
both with infinite combinations of unconditional and infinite
combinations of nonconglomerable conditional probabilities.
Coherence$_3$ leads to failures of the Dominance Principle only with
infinite combinations of nonconglomerable conditional
probabilities.\vadjust{\goodbreak}

A referee suggested that de Finetti might have been working with a
different Dominance Principle, here denoted Dominance*.

\begin{quote}
\emph{Dominance*}: Let $O_i$ and $O_j$ be two options in ${\cal O}$. If
$O_i$ uniformly (strictly) dominances $O_j$ and there exists an option
$O_k$ in ${\cal O}$ that is not itself dominated by some $O_t$ in
${\cal O}$, then $O_j$ is an inadmissible choice from ${\cal O}$.
\end{quote}
Dominance* requires that some option from ${\cal O}$ is
undominated if dominance signals inadmissibility. With respect to the
decision problems considered in this paper, each of our results
formulated with respect to the Dominance Principle obtains also with
Dominance*. Because Dominance* implies Dominance as we have
defined it, the only result that needs to be checked is
Theorem~\ref{thm:2}. In that case, so long as ${\cal O}$
contains options that correspond to a probability that satisfies the
law of total previsions in $\pi$ (as will all countably additive
probabilities) then Theorem~\ref{thm:3} says that such options will be
undominated. So, we could replace Dominance by Dominance* in the
results of this paper.

\begin{appendix}
\section{Finitely additive expectations}\label{app:lf}
This appendix gives the definitions of infinite prevision and finitely
additive expectation along with brief motivation for these definitions.
Details are given in
the supplemental article [Schervish, Seidenfeld and
Kadane (\citeyear{ssk2013})].
\subsection{Infinite previsions}\label{sec:infprev}
Our theorems assume that various random variables have finite
previsions. In the proof of Theorem~\ref{thm:1}, the possibility
arises that some other random variable has infinite prevision.
Definition~\ref{def:coh} makes no sense if infinite previsions are
possible. Fortunately, we can extend the concept of coherent$_1$
(conditional) prevision to handle infinite values, which correspond to
expressing a willingness either to buy or to sell a gamble, but not
both.
%
\begin{definition}\label{def:cohinf}
Let $\{P(X_i|B_i)\dvtx i\in I\}$ be a collection of conditional
previsions. The previsions are \emph{coherent$_1$} if, for every
finite $n$, every $\{i_1,\ldots,i_n\}\subseteq I$, all real
$\alpha_1,\ldots,\alpha_n$ such that $\alpha_j\leq0$ for all $j$
with $P(X_{i_j}|B_{i_j})=\infty$ and $\alpha_j\geq0$ for all $j$ with
$P(X_{i_j}|B_{i_j})=-\infty$, and all real $c_1,\ldots,c_n$ such that
$c_j=P(X_{i_j}|B_{i_j})$ for each $j$ such that $P(X_{i_j}|B_{i_j})$
is finite, we have
%
\begin{equation}
\label{eq:yacc} \inf_{\omega\in\Omega}\sum_{j=1}^n
\alpha_jB_{i_j}(\omega )\bigl[X_{i_j}(\omega
)-c_j\bigr]\leq0.
\end{equation}
That is, no linear combination of gambles may be uniformly strictly
dominated by the alternative option of abstaining.
\end{definition}
Notice the restrictions on the signs of coefficients in
Definition~\ref{def:cohinf}, namely that for each infinite prevision,
$\alpha_j$ has the opposite sign as the prevision. These restrictions
express the meaning of infinite
previsions as being \emph{one-sided} in the
sense that they merely specify that all real numbers are either
acceptable buy prices (for $\infty$ previsions) \emph{or} acceptable
sell prices
(for $-\infty$ previsions) but not fair prices for both transactions.
\citet{crisma97} and \citet{crisma01}
give alternate definitions of coherence for infinite previsions and
conditional previsions. But their definition does not make
clear the connection to gambling. However, the definition of Crisma, Gigante and
Millossovich (\citeyear{crisma97}) and Definition~\ref{def:cohinf} are equivalent for unconditional
previsions, as shown in the supplemental article [Schervish, Seidenfeld and
Kadane (\citeyear{ssk2013})].

\subsection{Prevision and expectation}\label{sec:prevexp}
Throughout this paper, an expectation with respect to a finitely
additive probability will be defined as a special type of linear
functional on a space of random variables. [See \citet{heathsudderth1978} for the case of bounded random variables.] Infinite
previsions are allowed in the sense of Section~\ref{sec:infprev}.
%
\begin{definition}\label{def:daniell}
Let $\L$ be a linear space of real-valued functions defined on
$\Omega$ that contains all constant functions, and
let $L$ be an extended-real-valued functional defined on $\L$.
If ($X,Y\in\L$ and $X\leq Y$) implies $L(X)\leq L(Y)$, we say that
$L$ is
\emph{nonnegative}.
We call $L$ an \emph{extended-linear functional} on $\L$, if, for all
real $\alpha,\beta$ and all $X,Y\in\L$,
%
\begin{equation}
\label{eq:el} L(\alpha X+\beta Y)=\alpha L(X)+\beta L(Y),
\end{equation}
whenever the arithmetic on the right-hand side of (\ref{eq:el}) is
well defined (i.e., not $\infty-\infty$) and where $0\times\pm
\infty=0$
in (\ref{eq:el}).
A nonnegative extended-linear functional is called a \emph{finitely
additive Daniell integral}. [See Schervish, Seidenfeld and
Kadane (\citeyear{ssk2008b}).]
If $L(1)=1$, we say that $L$ is \emph{normalized}. A~normalized
finitely additive Daniell integral is called a \emph{finitely
additive expectation}.
\end{definition}

Note that, if $\infty-\infty$ appears on the right-hand side of
(\ref{eq:el}), $L(\alpha X+\beta Y)$ still has a value, but the value
cannot be determined from (\ref{eq:el}).
Finitely additive expectations are essentially equivalent to
coherent$_1$ previsions, as we prove in the supplemental article.
Finitely additive expectations also behave like integrals in many
ways, as we explain in more detail in the supplemental article.
In particular, when the finitely additive expectation defined here is
restricted to bounded functions, it is the
same as the definition of integral developed by \citet{dunford-schwartz1958}, and it is the same as the integral used by \citet{dubins1975} in
his results about disintegrability. Hence,
Definition~\ref{def:daniell} is an extension of the definition of
integral from sets of bounded functions to arbitrary linear spaces of
functions.

\section{Proofs of results}\label{app:aux}
\subsection{Proof of Lemma \texorpdfstring{\protect\ref{lem:sp}}{1}}
Let $g$ be of the form (\ref{eq:gdef}). Let $P$ be such that $p=P(X)$
is finite, and let $q_0$ be such that $P[g(X,q_0)]$ is finite.
If $q\ne p$, then
%
\begin{equation}
\label{eq:pdiff} P\bigl[g(X,q)-g(X,p)\bigr]=\lambda\bigl((q,p)\bigr)\bigl[p-r(q,p,
\lambda)\bigr],
\end{equation}
according to (\ref{eq:diff}). Because Lebesgue measure is absolutely
continuous with respect to $\lambda$, neither $\lambda((q,p))$ nor
$p-r(q,p,\lambda)$ equals 0 and they have the same sign. It follows that
(\ref{eq:pdiff}) is
strictly positive. Since $p$ is finite, (\ref{eq:pdiff}) is finite
with $q=q_0$, so that $P[g(X,p)]$ is also finite and so $q=p$ provides
the unique minimum value of $P[g(X,q)]$.
\subsection{Proof of Lemma \texorpdfstring{\protect\ref{lem:ubf}}{2}}
Since $r(b,a,\lambda)=r(a,b,\lambda)$, it suffices to assume that
$a<b$. Let $\varepsilon>0$.

(i) If $\lambda_i((a,b))\geq\varepsilon$ and $b_0<b$ is such
that $\lambda_i((a,b_0))=\varepsilon$, then the probability obtained by
normalizing $\lambda_i$ on the interval $(a,b)$ stochastically dominates
the probability obtained by normalizing $\lambda_i$ on the interval
$(a,b_0)$. Hence, $r(a.b,\lambda_i)\geq r(a,b_0,\lambda_i)$. So, it
suffices to find a $\delta$ that implies
$r(a,b,\lambda_i)-a\geq\delta$ for all $i\in I$ and all $a<b$
such that $\lambda_i((a,b))=\varepsilon$. For the remainder of the
proof, let $a<b$ with $\lambda_i((a,b))=\varepsilon$, and let $Q$ be the
probability obtained by normalizing $\lambda_i$ on $(a,b)$.
Let $\lambda_0$ be $U$ times Lebesgue measure. Then
$\lambda_0((a,a+\varepsilon/U))=\varepsilon$, and
$r(a,a+\varepsilon/U,\lambda_0)=a+\varepsilon/(2U)$. Because $f_i\leq U$, it
follows that $Q$ stochastically dominates the probability obtained by
normalizing $\lambda_0$ on $(a,a+\varepsilon/U)$, hence $r(a,b,\lambda
_i)\geq
a+\varepsilon/(2U)$, and $r(a,b,\lambda_i)-a\geq\varepsilon/(2U)$. The
proof $b-r(a,b,\lambda_i)\geq\varepsilon/(2U)$ is similar, so
$\delta_\varepsilon$ can be taken equal to $\varepsilon/(2U)$.

(ii) Let $i\in I$, and assume that $\lambda_i((a,b))\geq\varepsilon$.
Since $f_j(v)>L_if_i(v)$ for all $v$, we have
$\lambda_j((a,b))\geq L_i\varepsilon$, so $\gamma_{i,\varepsilon}$ can be
taken to be
$L_i\varepsilon$.

\subsection{Proofs of Theorem \texorpdfstring{\protect\ref{thm:1}}{1} and
Corollary \texorpdfstring{\protect\ref{thm:1}}{1}}
Because a larger random variable has a larger prevision than a smaller
random variable, a necessary condition for (\ref{eq:dom}) is that
%
\begin{equation}
\label{eq:qsum} Z=P \Biggl[\sum_{i=1}^\infty
g_i(X_i,q_i) \Biggr]\leq P \Biggl[\sum
_{i=1}^\infty g_i(X_i,p_i)
\Biggr]<\infty.
\end{equation}
Hence, we will assume that $Z<\infty$
from now on. Also, it is necessary for (\ref{eq:dom}) that $q_i\ne
p_i$ for at least one $i$, so we will assume this also.

In light of (\ref{eq:diff}), we can write, for each finite $k>0$,
\begin{eqnarray*}
\infty>Z-W&=&P \Biggl[\sum_{i=1}^\infty
g_i(X_i,q_i)-\sum
_{i=1}^\infty g_i(X_i,p_i)
\Biggr]
\\
&=&\sum_{i=1}^k\lambda_i
\bigl((q_i,p_i)\bigr)[p_i-r_i]
+P \Biggl[\sum_{i=k+1}^\infty
g_i(X_i,q_i)-\sum
_{i=k+1}^\infty g_i(X_i,p_i)
\Biggr]
\\
&\geq&\sum_{i=1}^k\lambda_i
\bigl((q_i,p_i)\bigr)[p_i-r_i]-W,
\end{eqnarray*}
where the inequality follows because $g_i$ is nonnegative for each $i$
and where $r_i = r(q_i,p_i,\lambda_i)$ from \eqref{eq:ra2}.
Since $Z-W$ does not depend on $k$, it follows that
$\sum_{i=1}^\infty\lambda_i((q_i,p_i))[p_i-r_i]$ is finite.

Because of (\ref{eq2}) and (\ref{eq:qsum}), the two series
$\sum_{i=1}^\infty g_i(X_i(\omega),q_i)$ and\break  $\sum_{i=1}^\infty
g_i(X_i(\omega), p_i)$ are simultaneously finite with probability 1.
Let $B$ be the event that at least one of the two series is
finite. On $B^C$, both series sum to $\infty$, hence
(\ref{eq:dom}) fails unless $B^C=\varnothing$. Hence, we can assume
that $B=\Omega$ for the rest of the proof. It now follows that, for
all $\omega$,
%
\begin{eqnarray}
\label{eq:pq} &&\sum_{i=1}^\infty
g_i\bigl(X_i(\omega),q_i\bigr)-\sum
_{i=1}^\infty g_i
\bigl(X_i(\omega),p_i\bigr)
\nonumber
\\[-8pt]
\\[-8pt]
\nonumber
&&\qquad=\sum
_{i=1}^\infty\bigl[g_i\bigl(X_i(
\omega),q_i\bigr)-g_i\bigl(X_i(
\omega),p_i\bigr)\bigr].
\end{eqnarray}

We complete the proof by showing that
%
\begin{equation}
\label{eq5} P \Biggl[\sum_{i=1}^\infty
g_i(X_{i},q_{i})-\sum
_{i=1}^\infty g_i(X_{i},p_{i})
\Biggr]>0.
\end{equation}
Because a nonpositive random variable has nonpositive forecast,
(\ref{eq5}) implies that~(\ref{eq:dom}) cannot hold for all $\omega$.
In light of (\ref{eq:pq}), it suffices to show that
%
\begin{equation}
\label{eq:xyb} P \Biggl(\sum_{i=1}^\infty
\bigl[g_i(X_{i},q_{i})- g_i(X_{i},p_{i})
\bigr] \Biggr)>0.
\end{equation}

For each $k$,
%
\begin{equation}
\label{eq7} P \Biggl(\sum_{i=1}^k
\bigl[g_i(X_{i},q_{i})-g_i(X_{i},p_{i})
\bigr] \Biggr) = \sum_{i=1}^{k}
\lambda_i\bigl((q_i,p_i)\bigr)
(p_i-r_i) \geq0.
\end{equation}

Next, in light of (\ref{eq:diff}) and (\ref{eq7}), write
%
\begin{eqnarray}
\label{eq:split} &&P \Biggl(\sum_{i=1}^\infty
\bigl[g_i(X_i,q_i)-g_i(X_i,p_i)
\bigr] \Biggr)
\nonumber
\\[-8pt]
\\[-8pt]
\nonumber
&&\qquad= \sum_{i=1}^{k} \lambda_i
\bigl((q_i,p_i)\bigr) (p_i-r_i)
+
\nonumber
P \Biggl[\sum_{i=k+1}^{\infty}
\lambda_i\bigl((q_i,p_i)\bigr)
(X_i-r_i) \Biggr].
\end{eqnarray}
Since the left-hand side of (\ref{eq:split}) does not depend on $k$
and the
first sum on the right side is nondecreasing in $k$, it follows that
the second sum on the right-hand side is nonincreasing in $k$, and hence,
has a limit.
Let $T=\sum_{i=1}^\infty\lambda_i((q_i,p_i))(p_i-r_i)$, which is finite
and strictly positive (because $q_i\ne p_i$ for at least one $i$).
Then, the right-hand side of (\ref{eq:split}) becomes
%
\begin{equation}
\label{eq:sumT} T+\lim_{k\rightarrow\infty} P \Biggl[\sum
_{i=k+1}^\infty\lambda_i\bigl((q_i,p_i)
\bigr) (X_i-p_i) \Biggr].
\end{equation}
The proof will be complete if we can show that the limit in
(\ref{eq:sumT}) is 0.

First, we show that $\lim_{i\rightarrow\infty}\lambda_i((q_i,p_i))=0$.
If $\limsup_{i\rightarrow\infty}|\lambda_i((q_i,p_i))|= \ell>0$, then
there must exist a subsequence
$\{i_j\}_{j=1}^\infty$ with $|\lambda_{i_j}((q_{i_j},p_{i_j}))|>\ell
/2$ for
all $j$. For such a subsequence, the uniform spread condition implies
that there is $\delta_{\ell/2}>0$ such that $|p_{i_j}-r_{i_j}|\geq
\delta_{\ell/2}$. This would make $T=\infty$, a~contradiction.

It now follows that
\begin{eqnarray*}
\Biggl\llvert P \Biggl[\sum_{i=k+1}^\infty
\lambda _i\bigl((q_i,p_i)\bigr)
(X_i-p_i) \Biggr]\Biggr\rrvert &\leq& \max
_{i\geq k+1}\bigl|\lambda_i\bigl((q_i,p_i)
\bigr)\bigr|P \Biggl(\sum_{i=1}^\infty
|X_i-p_i| \Biggr)
\\
&=&V\max_{i\geq k+1}\bigl|\lambda_i\bigl((q_i,p_i)
\bigr)\bigr|,
\end{eqnarray*}
which can be made arbitrarily small by increasing $k$, and
(\ref{eq:xyb}) follows.

Corollary~\ref{thm:1} is equivalent to equation (\ref{eq5}), which is
established in the proof of Theorem~\ref{thm:1}.

\subsection{Proofs of Theorem \texorpdfstring{\protect\ref{thm:2}}{2} and
Corollary \texorpdfstring{\protect\ref{thm:2}}{2}}
Let $\pi= \{H_j\}_{j = 1}^\infty$ be a denumerable
partition. Nononglomerability means that there exists a random variable
$X$ such that either
\begin{eqnarray*}
\inf_j P(X|H_j) -P(X) &>& 0\quad\mbox{or}
\\
\sup_jP(X|H_j)-P(X)&<&0.
\end{eqnarray*}
Clearly, if $X$ satisfies one of
the above inequalities, $-X$ satisfies the other, hence we will
assume that the first inequality holds. Specifically, let
$p_X=P(X)$ and $p_j=P(X|H_j)$ for all $j$, and assume that
\[
\varepsilon= \inf_jp_j -p_X > 0.
\]
Also, for each $\omega\in\Omega$, let $i(\omega)$ be the unique
integer such
that $\omega\in H_{i(\omega)}$. Hence $H_j(\omega)=1$ if and only if
$j=i(\omega)$.

(2.1) Consider the following sum of individually fair options:
$X(\omega) - p_X$ and the countably many options
$-H_j(\omega)[X(\omega) - p_j]$ for $j=1,2,\ldots.$ Then, for each~$\omega$,
\begin{eqnarray*}
&&X(\omega) - p_X + \sum_{j=1}^\infty
-H_j(\omega)\bigl[X(\omega) - p_j\bigr]
\\
&&\qquad=X(\omega) - p_X - X(\omega) + p_{i(\omega)} =
-p_{X}+ p_{i(\omega)} \geq\varepsilon.
\end{eqnarray*}
Thus, the countable sum of the conditional forecasts for $X$ given
$H_j$, combined with the forecast for $X$ results in a
loss that is uniformly strictly dominated by 0.

(2.2) For an arbitrary set of forecasts $s_X$ for $X$ and $s_j$
for $X$ given
$H_j$ (for $j=1,\ldots$), the sum of the scores in state $\omega$ equals
%
\begin{eqnarray}
\label{eq16}
&&g_0\bigl(X(\omega),s_X
\bigr) + \sum_{j=1}^\infty H_j(
\omega )g_j\bigl(X(\omega ),s_j\bigr)\nonumber
\\
&&\qquad= g_0\bigl(X(\omega),s_X\bigr) +g_{i(\omega)}
\bigl(X(\omega), s_{i(\omega)}\bigr)
\\
&&\qquad=\int_{s_X}^{X(\omega)}\bigl[X(\omega)-v\bigr]\,d
\lambda_0(v) +\int_{s_{i(\omega)}}^{X(\omega)}\bigl[X(
\omega)-v\bigr]\,d\lambda_{i(\omega)}(v).
\nonumber
\end{eqnarray}

We can substitute the original forecasts $s_X=p_X$ and $s_j=p_j$,
$j=1,\ldots$
into~(\ref{eq16}) to obtain the total score
for each $\omega\in\Omega$. We can also identify dominating rival
forecasts $q_X$ and $q_j$, $j=1,\ldots,$ so that (\ref{eq16}) is
uniformly larger, for each state $\omega\in\Omega$ with $s_X=p_X$ and
$s_j=p_j$ than with $s_X=q_X$ and $s_j=q_j$.

Let $w_0=\lambda_0((p_X,p_X+\varepsilon))/2$, and let
$w_1=\gamma_{w_0}$, where $\gamma_{w_0}$ is from part (ii) of
Lemma~\ref{lem:ubf}. Let $q'$ be such that
$\lambda_0((q',p_X+\varepsilon))=w_0$. This makes
$\lambda_j((q',p_X+\varepsilon))\geq w_1$ for all $j$. For
each $j$, $p_j\geq p_X+\varepsilon$, so that $\lambda_j((q',p_j))\geq
w_1$. Let
$w_2=0.9\min\{w_0,w_1\}$, and let $q_j$ be such that
$\lambda_j((q_j,p_j))=w_2$ for all $j$. This makes $q_j> q'$ for all
$j$. Let $q_X$ be such that $\lambda_0((p_X,q_X))=w_2$. This makes $q_X<q'$.

We now form the difference between the scores for the original forecasts
and the rival forecasts. Subtracting (\ref{eq16}) with $s=q_X$ and
$s_j=q_j$ (for all $j$) from (\ref{eq16}) with $s=p_X$ and $s_j=p_j$
(for all $j$) yields
%
\begin{equation}
\label{eq17} \int_{p_X}^{q_X}\bigl[X(\omega)-v
\bigr]\,d\lambda_0(v)-\int_{q_{i(\omega
)}}^{p_{i(\omega)}}
\bigl[X(\omega)-v\bigr]\,d\lambda_{i(\omega)}(v).
\end{equation}
We need to find a positive number $\delta$ such that
(\ref{eq17}) is strictly greater than $\delta$ for all~$\omega$.
The difference in (\ref{eq17}) is greater than
\begin{eqnarray*}
&&\bigl[X(\omega)-q_X\bigr]\lambda_0
\bigl((p_X,q_X)\bigr)-\bigl[X(\omega)-q_{i(\omega)}
\bigr] \lambda_{i(\omega)}\bigl((q_{i(\omega)},p_{i(\omega)})\bigr)
\\
&&\qquad =w_2(q_{i(\omega)}-q_X)>w_2
\bigl(q'-q_X\bigr)>0.
\end{eqnarray*}
So, we set $\delta=w_2(q'-q_X)>0$, which
completes the proof.\vadjust{\goodbreak}

Corollary~\ref{thm:2} is immediate from {(2.2)} of
Theorem~\ref{thm:2}, as the existence of the rival set of dominating
forecasts, $\{q_j\}_{j=1}^\infty$, establishes that the forecaster
does not
minimize the infinite sum of expected scores by giving the forecast
$p_X$ and the conditional forecasts $\{p_j\}_{j=1}^\infty$.

\subsection{Proofs of Theorem \texorpdfstring{\protect\ref{thm:3}}{3} and
Corollary \texorpdfstring{\protect\ref{thm:3}}{3}}

(3.1) In order to show that (\ref{eq18}) cannot be uniformly strictly
positive, it is sufficient to show
%
\begin{equation}
\label{eq19} P \Biggl[\alpha_0(X-p_X)+\sum
_{j=1}^\infty\alpha_jH_j(X-p_j)
\Biggr]=0.
\end{equation}

Of course,
\[
P \Biggl[\alpha_0(X-p_X)+\sum
_{j=1}^\infty\alpha_jH_j(X-p_j)
\Biggr]= P\bigl[\alpha_0(X-p_X)\bigr]+P \Biggl[\sum
_{j=1}^\infty\alpha _jH_j(X-p_j)
\Biggr].
\]
Trivially,
%
\begin{equation}
\label{eq20} P\bigl[\alpha_0(X-p_X)\bigr]=0.
\end{equation}

Since $P$ satisfies the law of total previsions in $\pi$,
\[
P \Biggl[\sum_{j=1}^\infty
\alpha_jH_j(X-p_j) \Biggr]=P \Biggl[
P \Biggl[ \sum_{j=1}^\infty
\alpha_jH_j(X-p_i) \Big\rrvert \pi \Biggr]
\Biggr].
\]
For each $i$, $\sum_{j\ne i}\alpha_jH_j(\omega)[X(\omega)-p_j]=0$ for
all $\omega\in H_i$. It follows that, for every $i$,
\[
P \Biggl[\sum_{j=1}^\infty
\alpha_jH_j(X-p_j)\Big\rrvert H_i
\Biggr] =P\bigl[\alpha_iH_i(X-p_i)\mid
H_i\bigr],
\]
and trivially, $P[\alpha_iH_i(X-p_i)\mid H_i]=0$.

Thus,
\[
P \Biggl[\sum_{j=1}^\infty
\alpha_jH_j(X-p_i)\Big\rrvert \pi \Biggr]=0
\]
for all $\omega$, and it follows by the law of total previsions that
%
\begin{equation}
\label{eq23} P \Biggl[\sum_{j=1}^\infty
\alpha_jH_j(X-p_j) \Biggr] = 0.
\end{equation}
Equations (\ref{eq20}) and (\ref{eq23}) establish (\ref{eq19}).

(3.2) We must establish that there is no rival set of forecasts
$q_X$, and $\{q_j\}_{j=1}^\infty$ whose total score
uniformly dominates (\ref{eq24}). That is, there is no rival set of
forecasts such that for some $\varepsilon> 0$ and every $\omega$,
\begin{eqnarray*}
&&g_0\bigl(X(\omega),p_X\bigr)+\sum
_{j=1}^\infty H_j(\omega)g_j
\bigl(X(\omega),p_j\bigr) \\
&&\qquad\geq g_0\bigl(X(
\omega),q_X\bigr)+\sum_{j=1}^\infty
H_j(\omega)g_j\bigl(X(\omega),q_j\bigr) +
\varepsilon.
\end{eqnarray*}
It is sufficient to show that
%
\begin{eqnarray}
\label{eq26} \qquad &&P \Biggl\{g_0(X,q_X)+\sum
_{j=1}^\infty H_jg_j(X,q_j)
- \Biggl[g_0(X,p_X)+ \sum_{j=1}^\infty
H_jg_j(X,p_j) \Biggr] \Biggr\}
\nonumber
\\[-8pt]
\\[-8pt]
\nonumber
&&\qquad\geq0.
\end{eqnarray}
Write the left-hand side of (\ref{eq26}) as
%
\begin{equation}
\label{eq27} \qquad P\bigl[g_0(X,q_X)-g_0(X,p_X)
\bigr]+P \Biggl[\sum_{j=1}^\infty
H_jg_j(X,q_j)-\sum
_{j=1}^\infty H_jg_j(X,p_j)
\Biggr].
\end{equation}

That the first expectation in (\ref{eq27}) is nonnegative follows from
the fact that $g$ is strictly proper.
From the assumption that $P$ satisfies the law of total previsions in~$\pi$,
\begin{eqnarray*}
&& P \Biggl[\sum_{j=1}^\infty
H_jg_j(X,q_j)-\sum
_{j=1}^\infty H_jg_j(X,p_j)
\Biggr]\\
&&\qquad =P \Biggl[P \Biggl[\sum_{j=1}^\infty
H_jg_j(X,q_j)-\sum
_{j=1}^\infty H_jg_j(X,p_i)
\Big\rrvert \pi \Biggr] \Biggr].
\end{eqnarray*}
Using equation (\ref{eq:diff}) and the same logic as in part (3.1), we obtain, for each~$i$,
\begin{eqnarray*}
&& P \Biggl[\sum_{j=1}^\infty
H_jg_j(X,q_j)-\sum
_{j=1}^\infty H_jg_j(X,p_j)
\Big\rrvert H_i \Biggr]\\
&&\qquad=P\bigl[H_i \bigl\{
g_i(X,q_i)-g_i(X,p_i) \bigr
\}|H_i\bigr]
\\
&&\qquad=P\bigl[g_i(X,q_i)|H_i\bigr]-P
\bigl[g_i(X,p_i)|H_i\bigr]
\\
&&\qquad\geq 0,
\end{eqnarray*}
where the final inequality follows because $g_i$ is a proper scoring
rule and $P(\cdot|H_i)$ is a finitely additive expectation for all $i$.

Therefore, since $P$ satisfies the law of total previsions in $\pi$,
\[
P \Biggl[\sum_{j=1}^\infty
H_jg_j(X,q_i)-\sum
_{j=1}^\infty H_jg_j(X,p_i)
\Biggr]\geq0,
\]
which completes the proof of (\ref{eq27}).

Corollary~\ref{thm:3} is equivalent to the claim that for each set of
rival forecasts, $q_X$ and $\{q_j\}_{j=1}^\infty$, the second
prevision in (\ref{eq27}) is nonnegative, which was established
in the proof of (3.2).
\end{appendix}


\section*{Acknowledgments}
We thank Raphael
Stern, Department of Statistics, Carnegie Mellon University,
and two anonymous referees for
helpful advice about this paper.

\begin{supplement}[id=suppA]
\stitle{Infinite previsions and finitely additive expectations\\}
\slink[doi]{10.1214/14-AOS1203SUPP} 
\sdatatype{.pdf}
\sfilename{aos1203\_supp.pdf}
\sdescription{The expectation of a random variable $X$ defined on
$\Omega$ is
usually defined as the integral of $X$ over the set $\Omega$ with
respect to the underlying probability measure defined on subsets of
$\Omega$. In the countably additive setting, such integrals can be
defined (except for certain cases involving $\infty-\infty$)
uniquely from a probability measure on $\Omega$.
Dunford and Schwartz [(\citeyear{dunford-schwartz1958}), Chapter III] give a detailed analysis of
integration with respect to finitely additive measures that attempts to
replicate the uniqueness of integrals. Their analysis requires
additional assumptions if one wishes to integrate
unbounded random variables.
We choose the alternative of defining integrals as special types of
linear functionals. This is the approach
used in the study of the Daniell integral.
[See Royden (\citeyear{royden1968}), Chapter~13.] Then the measure of a set becomes the integral of its
indicator function. De Finetti's concept of
prevision turns out to be a finitely additive generalization of the
Daniell integral. (See Definition~\ref{def:daniell} in
Appendix \ref{sec:prevexp}.) 
We provide details on
the finitely additive Daniell integral along with details about the
meaning of infinite previsions and how to extend coherence$_1$\vadjust{\goodbreak} and
coherence$_3$ to deal with random variables having infinite
previsions. Infinite previsions invariably arise when dealing with
general sets of unbounded random variables.}
\end{supplement}

%



\printaddresses

\end{document}